\documentclass[12pt]{amsart}

\usepackage{amsmath}
\usepackage{amsfonts}
\usepackage{amssymb}
\usepackage{amsthm}
\usepackage{hyperref}
\usepackage{enumerate}
\usepackage{cleveref}
\usepackage{mathrsfs}
\usepackage[top=0.9in, bottom=1.15in, left=1.1in, right=1.1in]{geometry}
 \usepackage{hyperref}
\hypersetup{colorlinks=true,linkcolor=blue,citecolor=magenta}
\newtheorem{theorem}{Theorem}

\Crefname{conjecture}{Conjecture}{Conjectures}

\theoremstyle{plain}

\theoremstyle{plain}

\author{Aliyah Maxwell-Abrams and Robert Schneider}
\address{Department of Mathematical Sciences\newline
Michigan Technological University\newline
Houghton, Michigan 49931, U.S.A.}
\email{amaxwell@mtu.edu}
\email{robertsc@mtu.edu}

\title{Infinite series identities involving $\pi$ and $\operatorname{ln} 2$}
 
\begin{document}


\maketitle




\begin{center} {\it In celebration of Pi Day 2025}\end{center}





\section{Introduction}

 In this paper we prove identities for six infinite series whose   values involve linear combinations of $\pi$ and $\operatorname{ln} 2$, that do not appear in standard infinite series references.

\section{New formulas for $\pi$}
\label{sec:goodstory}

\noindent In a  previous paper \cite{pi/3}   published to celebrate Pi Day 2022, the second author   proves the identity
\begin{equation}\label{OG}
\sum_{n=1}^{\infty} \frac{1}{n(2n-1)(4n-3)}\  =\  \frac{\pi}{3},\end{equation}
as well as a list of twenty-seven further  infinite series identities involving linear combinations of $\pi, \pi^2,\operatorname{ln} 2$ and $\operatorname{ln} 3$. 

In the closing of \cite{pi/3}, one more formula is noted as a  problem for the reader, which is to find the value of a similar series absent from the list:
\begin{equation}\label{question}
\sum_{n=1}^{\infty} \frac{1}{n(4n-1)(4n-3)}.\end{equation}
In a   2024 undergraduate research project at Michigan Technological University mentored by the second author, the first author of this paper found the  value of series \eqref{question}.

\begin{theorem}\label{thm1} We have the identity
\begin{equation*}
\sum_{n=1}^{\infty} \frac{1}{n(4n-1)(4n-3)}\  =\  \frac{\pi}{3}-\operatorname{ln} 2. \end{equation*}
\end{theorem}

To prove this theorem, we combine entries (15) and (16) of \cite{pi/3}. Let $S_{15}$ and $S_{16}$ denote the infinite series given by these two entries, respectively:\footnote{The identity for $S_{15}$ is proved in \cite{figurate} as the decagonal case of a  formula related to figurate numbers.}

\begin{equation}\label{seeds1}
S_{15}=\sum_{n=1}^{\infty} \frac{1}{n(4n-3)}=\frac{\pi+6\operatorname{ln}2}{6},\  \  \  \  \  \   S_{16}=\sum_{n=1}^{\infty} \frac{1}{n(4n-1)}=\frac{6\operatorname{ln}2 - \pi}{2}. 
\end{equation}
These series are (absolutely) convergent by either the integral or series comparison test. 
We note     these identities were proved by the second author using combinations of the following (conditionally convergent) series, carefully treating the series as limits of combinations of Maclaurin series:\footnote{They are the limiting cases of $\operatorname{ln}(1+x)$ as  $x\to 1^-,$ and $\operatorname{arctan} x$ as $x\to 1,$ respectively. The reader is referred to    \cite{history, vid} for details about  Taylor series representations for $\pi$.}
\begin{flalign}\label{classical}
1-\frac{1}{2}+\frac{1}{3}-\frac{1}{4}+\dots\  =\  \operatorname{ln} 2,\  \  \  \  \  \  
1-\frac{1}{3}+\frac{1}{5}-\frac{1}{7}+\dots\  =\  \frac{\pi}{4}.
    \end{flalign}
The second series above is known as the M\={a}dhava-Gregory-Leibniz series. 

To prove the theorem, we  compute $\frac{1}{2}(S_{15}-S_{16}),$ with      $S_{15}$ and $S_{16}$ as in \eqref{seeds1} above. We note that $S_{15}-S_{16}$ is equal to
\begin{flalign}
\sum_{n=1}^{\infty} \frac{1}{n(4n-3)}-\sum_{n=1}^{\infty} \frac{1}{n(4n-1)}= \sum_{n=1}^{\infty} \frac{(4n-1)-(4n-3)}{n(4n-1)(4n-3)} =2\sum_{n=1}^{\infty} \frac{1}{n(4n-1)(4n-3)}.
\end{flalign}
Dividing by two, with a little arithmetic on the right-hand side, completes the proof. 

With a little arithmetic, the identities in \eqref{OG} and Theorem \ref{thm1} give   new formulas for $\pi$:
\begin{flalign}\label{OG2}
\pi\  &=\  3\sum_{n=1}^{\infty}\frac{1}{n(2n-1)(4n-3)}\  =\  3\ln 2 + 3\sum_{n=1}^{\infty} \frac{1}{n(4n-1)(4n-3)}.
\end{flalign}

\section{Further formulas}

With the identity in Theorem \ref{thm1} now in hand, we set our sights on supplementing the list of identities given in  \cite{pi/3} with further identities of similar forms. Theorem \ref{thm1} is a low-hanging fruit that was directly asked for in the previous paper. Upon further inspection, we were able to spot other gaps in the list of twenty-seven identities. 

The infinite series listed in  \cite{pi/3} are classified according to the number of linear factors in the denominators of the summands. Observing there were certain factors like $4n+1$ that showed up only rarely in the denominators of series in the  list of identities in \cite{pi/3}, we proved new identities of similar types. 

\begin{theorem}\label{thm2} We have the following identities: 
\begin{enumerate}[(a)] 
\item $$\sum_{n=1}^{\infty} \frac{1}{n(2n-1)(4n-1)(4n-3)}\  =\ 2\operatorname{ln} 2-\frac{\pi}{3},$$
\item $$\sum_{n=1}^{\infty} \frac{1}{n(4n+1)(4n-1)(4n-3)}\  =\  \frac{\pi}{12}-\operatorname{ln} 2+\frac{1}{2},$$
\item $$\sum_{n=1}^{\infty} \frac{1}{(2n-1)(4n+1)(4n-1)(4n-3)}\  =\ \frac{1}{6}\left(\operatorname{ln}2-\frac{\pi}{4}+\frac{1}{2}\right),$$
\item $$\sum_{n=1}^{\infty} \frac{1}{n(2n-1)(4n+1)(4n-3)}\  =\  \frac{2}{3}\left(1-\operatorname{ln} 2\right),$$
\item $$\sum_{n=1}^{\infty} \frac{1}{n(2n-1)(4n+1)(4n-1)(4n-3)}\  = \frac{1}{3}\left(4\operatorname{ln}2-\frac{\pi}{2}-1\right).$$
\end{enumerate}
\end{theorem}

\


These identities are proved   similarly to Theorem \ref{thm1} above. Let $S_1$ denote the infinite series on the left-hand side of Theorem \ref{thm1}. 
Let $S_{{a}}, S_{{b}}, S_{{c}}$, etc., denote the series in identities ({\it a}), ({\it b}), ({\it c}), etc., of Theorem \ref{thm2} above, respectively. Much like $S_{15}, S_{16}$ defined in \eqref{seeds1},  let $S_{17}, S_{19},$ denote the series in identities (17) and (19) of \cite{pi/3}: 
\begin{flalign}\label{seeds2}
S_{17}&=\sum_{n=1}^{\infty} \frac{1}{(2n-1)(4n-1)(4n-3)}=\frac{\operatorname{ln}2}{2},\\  S_{19}&=\sum_{n=1}^{\infty} \frac{1}{(4n+1)(4n-1)(4n-3)}=\frac{\pi-2}{16}.
\end{flalign}
These are proved by  linear  combinations of   identities   \eqref{OG}, \eqref{seeds1} and \eqref{classical}. By similar steps to the proof of Theorem \ref{thm1}, namely,    combining series term-wise using fraction arithmetic,  the    identities in Theorem \ref{thm2} follow from these formulas:
 \begin{flalign} 
S_{{a}}&= 2S_{17}-S_1, \\
S_{{b}}&=S_1-4S_{19},\\
S_{c}&=\frac{1}{3}(S_{17}-2S_{19}),\\
S_{d}&= 8S_{c}-S_{a},\\
%
S_{{e}}&=2S_{c}-S_{b}.
\end{flalign}
We note $S_{d}$ and $S_{e}$ both depend on   preceding identities $S_{a}$ and $S_{c}$ in Theorem \ref{thm2}. As an example of one of these proofs, let us derive $S_{d}$ explicitly. Note that $8S_{c}-S_{a}$ is equal to  
\begin{flalign}
\sum_{n=1}^{\infty}\frac{8n}{n(2n-1)(4n+1)(4n-1)(4n-3)}&-\sum_{n=1}^{\infty}\frac{4n+1}{n(2n-1)(4n+1)(4n-1)(4n-3)}\\
\nonumber &=\sum_{n=1}^{\infty}\frac{4n-1}{n(2n-1)(4n+1)(4n-1)(4n-3)}=S_{d}.
\end{flalign}
The same arithmetic steps, applied to the right-hand sides of ({\it a}) and ({\it c}), lead to  ({\it d}).

\begin{section}{Conclusion}
    
The identities we prove above do not appear in standard tables of infinite series at our disposal \cite{tables, Jolley}, nor in popular books about $\pi$ such as \cite{history,   joy}, nor among the well-known infinite series identities of Euler   (see e.g. \cite{Dunham, Maor}). 

We note that identities such as these will not give approximations to $\pi$ that converge rapidly enough to be useful for computations, by comparison to   formulas used to compute billions upon billions  of digits of $\pi$ (see \cite{Wikipedia1}), but they do converge reasonably quickly -- generally more so in a given series, the greater the number of factors  in the denominator. For example,    $S_{19}$ can be used to compute $\pi$:
\begin{equation}\label{approx}
\pi \  =\  2+16S_{19}. 
\end{equation}
Summing just the first nine terms of $S_{19}$, equation \eqref{approx} yields an approximation $\pi \approx 3.140133\dots$ accurate to two decimal places, and the first 217 terms of $S_{19}$ gives $\pi \approx 3.141590\dots$ accurate to five   places. By contrast, to use the M\={a}dhava-Gregory-Leibniz series $\pi=4\cdot (1-1/3+1/5-1/7+\dots)$ in \eqref{classical}   to produce an   estimate accurate to five decimal places takes  $500{,}000$ terms \cite{Wikipedia}. On the other hand,  Ramanujan proved an efficient formula whose $n=1$ summand alone gives the first seven digits of $\pi$, with  each additional term of the summation giving eight additional digits of the decimal expansion \cite{Wikipedia2}.\footnote{For     milestones in computing digits of $\pi$,   see \cite{Wikipedia2}.}

We encourage the  reader to seek further identities to supplement the lists of infinite series formulas made in \cite{pi/3} and in this paper, perhaps combining these formulas with known identities from the literature. Excellent starting points would be Dunham's survey \cite{Dunham} of the infinite series methods of Euler, and Beckmann's historical survey  \cite{history}  of  formulas for $\pi$. For further reading, see \cite{figurate} in which a theory is developed of identities similar to those in \cite{pi/3} and in this work. 
\end{section}

\section*{Acknowledgements}

This paper was written as part of an undergraduate research project undertaken by the first author at Michigan Technological University, Houghton, Michigan, U.S.A., mentored by the second author (Fall 2024). The authors are   grateful to the MTU McNair Scholars Program, MTU's Summer Undergraduate Research Fellowship (SURF), and the MTU Department of Mathematical Sciences, for supporting our joint research projects (2022 -- 2024). Thank you to James Sellers for advice that improved this paper.

\end{document}